\documentclass[%
  a4paper,
  onecolumn,
  colorlinks,
]{preprint}

\usepackage{silence}
\usepackage[english]{babel}

\usepackage[noadjust]{cite}

\usepackage{graphicx}
\usepackage{pgf}

\usepackage{amsfonts}
\usepackage{amsmath}
\usepackage{amssymb}
\usepackage{amsthm}
\usepackage{makecell}

\usepackage{todonotes}

\theoremstyle{plain}
\theoremstyle{plain}\newtheorem{lemma}{Lemma}
\theoremstyle{definition}\newtheorem{remark}{Remark}

\DeclareMathOperator\sftmx{softmax}
\DeclareMathOperator\Fenc{F_e}
\DeclareMathOperator\Fclstr{F_c}

\newcommand{\R}{\ensuremath{\mathbb{R}}}
\newcommand{\N}{\ensuremath{\mathbb{N}}}

\newcommand{\tA}{\ensuremath{\widetilde{A}}}

\newcommand{\tx}{\ensuremath{\tilde{x}}}

\def\conv{\ensuremath{\mathsf{Conv}}}
\def\elu{\ensuremath{\mathsf{ELU}}}
\def\fc{\ensuremath{\mathsf{FC}}}
\def\tmat{\ensuremath{\mathsf{T}}}
\newcommand{\te}{\ensuremath{t_{\operatorname{e}}}}

\newcommand{\xs}{\ensuremath{x^{\ast}}}

\newcommand{\astr}{{\alpha^*}}
\newcommand{\astro}{{\beta^*}}
\newcommand{\astrt}{{\delta^*}}
\newcommand\POD[1]{\texttt{POD-#1}}
\newcommand\PAEqr[2]{\texttt{PAE-#1.#2}}
\newcommand\PAEqrp[3]{\texttt{PAE-#1.#2(p=#3)}}

\newcommand{\trans}{\ensuremath{\mkern-1.5mu\mathsf{T}}}

\def\nclstr{q}
\def\xrdr{p}

\begin{document}

\title{Deep polytopic autoencoders for low-dimensional linear parameter-varying
  approximations and nonlinear feedback controller design}

\author[$\ast$]{Jan Heiland}
\affil[$\ast$]{Department of Mathematics and Natural Sciences,
  Institute of Mathematics, Technische Universit{\"a}t Ilmenau,
  Germany.\authorcr
  \email{jan.heiland@tu-ilmenau.de}, \orcid{0000-0003-0228-8522}}

\author[$\dagger$]{Yongho Kim}
\affil[$\dagger$]{Max Planck Institute for Dynamics of Complex Technical
  Systems, Sandtorstra{\ss}e 1, 39106 Magdeburg, Germany.
  \email{ykim@mpi-magdeburg.mpg.de}, \orcid{0000-0003-4181-7968}}

\author[$\ddagger$]{Steffen W. R. Werner}
\affil[$\ddagger$]{Department of Mathematics and
  Division of Computational Modeling and Data Analytics,
  Academy of Data Science, Virginia Tech,
  Blacksburg, VA 24061, USA.\authorcr
  \email{steffen.werner@vt.edu}, \orcid{0000-0003-1667-4862}}

\shorttitle{Polytopic autoencoders for LPV-based nonlinear control}
\shortauthor{Heiland, Kim, Werner}
\shortdate{2025-01-23}

\abstract{%
  Polytopic autoencoders provide low-di\-men\-sion\-al parametrizations of states
in a polytope.
For nonlinear PDEs, this is readily applied to low-dimensional linear
parameter-varying (LPV) approximations as they have been exploited for efficient
nonlinear controller design via series expansions of the solution to the
state-dependent Riccati equation.
In this work, we develop a polytopic autoencoder for control applications and
show how it improves on standard linear approaches in view of LPV
approximations of nonlinear systems.
We discuss how the particular architecture enables exact representation of
target states and higher order series expansions of the nonlinear feedback law
at little extra computational effort in the online phase and how the linear
though high-dimensional and nonstandard Lyapunov equations are efficiently
computed during the offline phase.
In a numerical study, we illustrate the procedure and how this approach can
reliably outperform the standard linear-quadratic regulator design.

}

\novelty{}



\maketitle



\section{Introduction}%
\label{sec:intro}

The design of feedback controllers for nonlinear high-dimensional systems is a
challenging task and there exists no generally applicable and computationally
feasible approach.
Here, we consider input-affine systems of the form
\begin{equation} \label{eq:gen-affine-system}
  \dot x(t) = f(x(t)) + Bu(t), \quad y(t) = C x(t),
\end{equation}
where for time $t \geq 0$, $x(t)\in \R^{n}$ denotes the state, $u(t) \in \R^{m}$
and $y(t)\in \R^{\ell}$ denote the system's input and output,
$f\colon \R^{n}\to \R^{n}$ is a possibly nonlinear function,
and $B \in \R^{n \times m}$ and $C \in \R^{p \times n}$ are linear input and
output operators represented as matrices.
The general idea of feedback controller design is to identify a mapping
$\mathcal{K}\colon x(t) \to u(t)$ so that with the choice of
$u = \mathcal{K}(x)$, the system's state $x$ and output $y$ satisfy desired
criteria.
A typical application is feedback stabilization that seeks to achieve
$x(t) \to 0$ as $t \to \infty$ or the more general task of
\emph{set-point control}, namely $x(t) \to \xs$ for some state
$\xs \in \R^{n}$ that fulfills~$f(\xs) = 0$.

For the controller design for systems of the form~\cref{eq:gen-affine-system},
methods like backstepping~\cite{Kok92},
feedback linearization~\cite[Chap.~5.3]{Son98}, and
sliding mode control~\cite{Dod15} require structural assumptions,
whereas general approaches based on the Hamilton-Jacobi-Bellman (HJB) equations
are only feasible for very moderate system sizes or in combination with model
reduction; see, e.g.,~\cite{BreKP19} for
polynomial expansions,~\cite{DolKK21} for tensorized approximations,
or~\cite{AllS20} for an approach that addresses general PDEs based on
first reducing the state-space dimension through, e.g., POD.
Other efforts try to approximate the solution of the HJB equations by deep
neural networks; see~\cite{DarLM20}.
The commonly employed state-dependent Riccati equation (SDRE) approximation to
the HJB equations (see, e.g.,~\cite{Cim12} for a survey) lifts a significant
portion of the computational complexity.
Nonetheless, this approach became available for large-scale systems only
recently due to further simplifications; see~\cite{BenH18} for a linear
update scheme and~\cite{AllKS23} for approximations through series expansions
for systems of a certain parametric structure.
These structural assumptions have been lifted in our previous work
(see~\cite{morHeiW23}) so that, with the help of model reduction, now a large
class of nonlinear control systems can be treated with series expansions of
the SDRE.

The underlying idea is the approximation of the nonlinear
system~\cref{eq:gen-affine-system} through a very low-dimensional linear
parameter-varying (LPV) system so that the series expansion of any nonlinear
feedback law can be formulated with respect to a few parameters rather than
the full state;
see our previous works~\cite{morHeiW23, morDasH23} and related
discussions in~\cite{AllKS23, AllP24}.
In particular, the results in~\cite{morHeiW23} indicate a measurable
improvement of the nonlinear SDRE feedback over a related, well-tuned and
robust linear-quadratic regulator (LQR) but left space for further performance
improvements, for example, through
\begin{itemize}
  \item[(i)] higher order expansions of the nonlinear feedback, or
  \item[(ii)] nonlinear parametrizations of the velocity states.
\end{itemize}
For higher order expansions of nonlinear feedback laws, a low-dimensional
parametrization of the states is crucial since the number of coefficients
grows drastically with the dimension of the parametrization and the order of
the expansion so that both topics are linked to the single question if a
nonlinear parametrization can provide sufficiently accurate low-dimensional
approximations with only few degrees of freedom.

It is known that linear parametrizations like the (linearly) optimal approach
of \emph{proper orthogonal decomposition} (POD, see, e.g.,~\cite{FreM22}) have a
natural limit on the approximation accuracy, which is commonly expressed
through the Kolmogorov $n$-width; see~\cite{OhlR13}.
Thus, there have been many efforts and positive results concerning the use
of nonlinear model order reduction approaches including, for example,
polynomial manifold maps~\cite{BarF22, BucGH23} or general nonlinear
autoencoders~\cite{FreM22, SheYSetal21, RizAV18}.
An alternative for accurate low-dimensional parametrizations is found in
the use of local bases~\cite{AmsZF12}; see the recent works that consider local POD
bases state reconstruction on subsets of the given
dataset \cite{AmsH16, AmZa12, KaSh22} or our previous work on local
bases approximations via autoencoders \cite{HeiK24,morHeiKim24}.
Basically, local bases are designed for subsets of states considered to have
high similarities.
However, local bases require a mechanism to identify the subsets and which basis
is to be used.
Standard approaches use clustering for this purpose, for example, $k$-means
clustering, fuzzy $c$-means clustering, and deep
clustering~\cite{ArVa07, Du73, FaTh20}.

In this work, we call on autoencoders with smooth basis selection and present a
purely computational approach to design reduced-order parametrizations of
nonlinear equations and provide algorithms and basic results for applications
in nonlinear feedback control.
Therefore, we extend our recently developed framework of polytopic
autoencoders~\cite{HeiK24} that combines an autoencoder with the selection
of local bases in terms of polytope vertices towards the use in set-point
control.
In particular, we directly impose the set-point condition $f(\xs)=0$ onto the
autoencoder so that we obtain a nonlinear encoding of very small dimension with
arbitrary accuracy in the relevant part of the state space.
To this end, we derive explicit expressions for the Jacobian of the decoder,
which we then can exploit in theory and practice for analysis and synthesis of
higher order expansions of the SDRE feedback.
The precomputation of the expansion coefficients requires the
solution of large-scale matrix-valued algebraic equations of a nonstandard form.
Here, we call and expand on recent advances in the solution of
indefinite matrix Lyapunov and Riccati equations; see~\cite{LanMS14, SaaW24}.
We note that the derived representation through a polytopic LPV system is
suited for controller design on the base of linear matrix inequalities (LMIs);
see, e.g.,~\cite{ApkGB95}.
However, since the reduction only considers the structure but not the
high-dimensional state space, the systems to solve are far too large for
state-of-the-art LMI solvers.
We exemplify the procedure and potentials of this general-purpose approach to
nonlinear feedback design with a numerical example based on the two-dimensional
flow past a cylinder.

The overall structure of this paper can be summarized as follows:
We recall in \Cref{sec:ld-lpv} that a general nonlinear control system
can be approximated by a low-dimensional LPV system, which can be
exploited for controller design in particular if the parameter dependency
itself is linear.
In view of providing a performant but maximally low-dimensional
parametrization, we resort to clustered polytopic autoencoders that encode
the states of a dynamical system in a polytope and decode them with locally
adapted bases; see \Cref{sec:paes}.
An interesting result in view of controller design is that this nonlinear
decoding has a Taylor expansion with vanishing second-order terms
(\Cref{lem:taylor-of-crhorho}).
For the actual controller design in \Cref{sec:hoSDRE-expansion}, we
consider series expansions of state-dependent Riccati equation.
Here, we extend the formulas to consider nonlinear parameter dependencies and
higher order terms in the expansion.
We discuss and showcase, how the relevant large-scale matrix equations can be
solved numerically.
Finally, we report numerical experiments that illustrate the performance
improvement of the polytopic autoencoder over standard approximation schemes
that explore the different levels of approximation to the underlying nonlinear
SDRE feedback design, and that investigate the gain in controller performance
for a commonly considered flow control problem; see \Cref{sec:num-exps}.
As a result, we obtain a model-based nonlinear controller of reduced complexity,
with the reduction tailored to snapshot data.
The paper is concluded in \Cref{sec:conclusions}.


\section{Low-dimensional linear parameter-varying approximations}%
\label{sec:ld-lpv}

For nonlinear systems like~\cref{eq:gen-affine-system}, we consider
\emph{quasi}-LPV approximations of the form
\begin{equation} \label{eq:nonlsys-lpv}
  \dot{x}(t) = A\bigl( \rho(x(t)) \bigr) x(t) + Bu(t), \quad y(t) = Cx(t).
\end{equation}
Such systems can be derived from an (exact) state-dependent coefficient
representation such as $f(x) = \tA(x) x$ and a parametrization of the state
$x(t) \approx \tilde x(t) = \nu(\mu(x(t))$, where $\mu$ is an encoder and
$\nu$ is a decoder; see~\cite{morHeiW23}.
In this case, we set $\rho = \mu(x)$ and $A(\rho) := \tA(\nu(\rho))$.
Note that the parameter $\rho$ depends on the state,
hence the term \emph{quasi}-LPV.

\begin{remark}\label{rem:Alinear}
  In the example application of fluid flow models, the coefficient function
  $\tA$ is readily available; cf.~\cite{HeiBB22}.
  Here, the function $\tA$ is linear so that a linear decoding will lead to
  linear parameter dependencies as they have been proven beneficial for
  controller design; see, e.g.,~\cite{ApkGB95}.
\end{remark}

In what follows, we refer by $\rho_{i}$ to the $i$-th component of a
vector-valued quantity $\rho$.
To reduce the notation, we omit dependencies and sometimes write, e.g.,
$\rho$ or $\rho(t)$ instead of $\rho(x(t))$.
Furthermore, we will denote the \emph{Kronecker product} by $\otimes$.


\section{Polytopic parametrizations for control}%
\label{sec:paes}

In polytopic autoencoders (PAEs), the state is represented through generalized
barycentric coordinates with respect to $R \in \N$ vertices
$\{ \tilde{x}_{k} \}_{k = 1}^{R} \subset \R^{n}$ of a polytope;
see, e.g., our previous work~\cite{HeiK24} for a more detailed introduction.
As a polytope is bounded, one can always assume without loss of generality that
the coordinates $\rho(t) \in \R^{R}$ in the reconstruction
will be positive and sum up to one, i.e., for the $\tilde x(t)$ in the polytope
it holds that
\begin{equation} \label{eq:polytopic-reconstruction}
  \tilde{x}(t) = \sum_{k = 1}^{R} \rho_{k}(t) \tilde{x}_{k},
  \quad\text{with } \rho(t) \geq 0 \text{ and }
    \sum_{k = 1}^{R} \rho_{k}(t) = 1.
\end{equation}
An immediate advantage of the polytopic parametrization for the intended
controller design through series expansions results from the boundedness of the
reconstruction, which keeps the controller gains bounded so that we can expect
good approximations in series expansions while possible scaling issues will be
detected when precomputing the expansion coefficients.

For the application in controller design, the consistency at the target
state $\xs$ is important.
We note that a design like~\cref{eq:polytopic-reconstruction} is
not capable of exactly parametrizing a target value $\xs = 0$ unless the
vertices $\tx_{k}$ come with a specific linear dependency.
Moreover, for the intended series expansions around the target state, accuracy
is guaranteed if $\|\rho(t)\|$ is small, which, however, contradicts the
summation condition in \eqref{eq:polytopic-reconstruction}.
To address these two requirements, we include $\tx_{0} = 0$ in the training of
the reconstruction basis, where we impose the positivity and summation
constraint on the $\rho$ variable extended by $\rho_{0}$, i.e.,
\begin{equation} \label{eq:theta-positive-sumone}
  \rho_{0}(t) + \sum_{k = 1}^{R} \rho_{k}(t) = 1.
\end{equation}
This way, the tuple $(1, \rho(t)) = (1, 0) \in \R^{R + 1}$ will be a
feasible value in the extended parameter space that exactly represents the
setpoint $\xs = 0 \in \R^{n}$.
Moreover, smoothness in the parameter space, will ensure that $\rho$ is small
around the origin.

As for the use of local bases, we assume that $q$ different polytopes with vertices
\begin{equation*}
  \{ \tilde{x}_{k}^{(j)}\}_{k = 1, \dotsc, R; \, j=1,\dotsc, q} \subset \R^{n}
\end{equation*}
are to be used as local bases for the reconstruction.
Then, a direct implementation would compute both the local coordinates $\rho$
and the index $j$ of the relevant polytope and reconstruct according
to the formula
\begin{equation}\label{eq:clustered-polytopic-reconstruction}
  \tilde x = \begin{bmatrix}
    \tilde x_1^{(1)} & \hdots & \tilde x_R^{(1)} &
    \tilde x_1^{(2)} & \hdots & \tilde x_R^{(2)} &
    \hdots &
    \tilde x_1^{(q)} & \hdots & \tilde x_R^{(q)}
  \end{bmatrix}
  (\alpha \otimes \rho)
\end{equation}
where, in this case, $\alpha = e_{j}$ is the $j$-th canonical basis vector.
Instead of this \emph{hard selection} of the $j$-th polytope, we rather use a
smooth clustering network (e.g.,~\cite{YanLL12GM, FaTh20}) that blends between
the polytopes by choosing $\alpha = c(\rho) \in \R^{q}$, with
$c(\rho) \geq 0$ and $\sum_{\ell = 1}^{q} c(\rho)_{\ell} = 1$.
We note that the summation condition on $c(\rho)$ allows to consider
$c(\rho) \otimes \rho$ as generalized barycentric coordinates in
the polytope spanned by all vertices as
in~\cref{eq:clustered-polytopic-reconstruction}; see~\cite[Lem~3.1]{HeiK24}.

Under the preceding considerations, we propose a polytopic
autoencoder designed for set-point control towards $\xs = 0$ that consists of
the following three components:
\begin{enumerate}
  \item[(i)] A nonlinear encoder $\mu\colon \R^{n} \to \R^{r + 1}$ so that
    \begin{equation*}
      (\rho_{0}, \rho) = \mu(x) =: \sftmx(\Fenc(x) + e_{1}),
    \end{equation*}
    with a scalar function $\rho_{0}$ and $\rho(t) \in \R^{r}$, where
    $e_{1} \in \R^{r+1}$ is the first unit vector and
    $\Fenc\colon \R^{n} \to \R^{r+1}$ is a nonlinear
    convolutional neural network described via
    \begin{align*}
      h_{1} & = \tmat x,\\
      h_{2} & = \elu(\conv_{1}(h_{1})),\\
      & \hspace{.5em}\vdots \\
      h_{L} & = \elu(\conv_{L-1}(h_{L-1})),\\
       \Fenc(x) & = \fc(h_{L}).
    \end{align*}
    Here, $\tmat$ denotes the interpolation of spatially distributed data onto a
    pixel grid as it is needed for the convolutions,
    $\elu$ is the exponential linear unit (ELU) activation
    function~\cite{CleUH15},
    $\fc$ is a fully connected layer, and
    $\conv_{i}$ is a convolutional layer, $i\in\{1,2,\cdots , L-1\}$ for some
    depth parameter $L\in \N$.
    Each layer is designed without bias terms to ensure $\Fenc(0) = 0$
    so that with the \emph{softmax function}
    \begin{equation*}
      \sftmx(\xi)_{i} = \tfrac{\xi_{i} \tanh(10 \xi_{i})}{\sum_{j = 1}^{r+1}
        \xi_{j} \tanh(10 \xi_{j})}
    \end{equation*}
    for $i = 1, \ldots , r+1$ and
    $\xi=(\xi_{1}, \ldots , \xi_{r+1}) \in \R^{r+1}$,
    we achieve, by design, that the output $(\rho_{0}(t), \rho(t))$ satisfies
    \begin{equation*}
      (\rho_{0}(t), \rho_{k}(t)) \geq 0 \quad\text{and}\quad
      \rho_{0}(t)+ \sum_{k = 1}^{r} \rho_{k}(t) = 1;
    \end{equation*}
    and in particular $\mu(0) = (1, 0)$.

  \item[(ii)] A smooth clustering network $c\colon \R^{r} \to \R^{\nclstr}$,
    with
    \begin{equation} \label{eq:c-clustering-net}
      c(\rho) = \sftmx (\Fclstr(\rho) + \hat{e}_{1}),
    \end{equation}
    where $\hat{e}_{1}$ is the first unit vector in $\R^{\nclstr}$ and
    \begin{equation*}
      \Fclstr(\rho) = \elu(\fc(\rho)),
    \end{equation*}
    with no bias, i.e., $\Fclstr(0) = 0$.
    By design, for the output $c(\rho)$, we have that
    \begin{equation*}
      c(\rho(t)) \geq 0
      \quad\text{and}\quad
      \sum_{i = 1}^{\nclstr} c(\rho(t))_{i} = 1
    \end{equation*}
    hold, so that $(\rho_{0}, c(\rho) \otimes \rho)$
    satisfies the positivity and summation
    conditions~\cref{eq:theta-positive-sumone}.

  \item[(iii)] A trainable linear map $W \in \R^{n \times \nclstr r}$ that
    realizes the reconstruction
    (cf.~\cref{eq:clustered-polytopic-reconstruction}) as
    \begin{equation} \label{eq:reconstruction-w-W}
      \tx = W \big( c(\rho) \otimes \rho \big).
    \end{equation}
\end{enumerate}

A few comments on the design of the autoencoder and motivation:
\begin{itemize}
  \item Generally, the property that $\Fenc(0) = 0$ and $\Fclstr(0) = 0$ can
    be achieved by using standard (convolutional) linear layers without
    bias terms and with any activation function that crosses the origin.

  \item The use of the extra $\rho_{0}$ variable ensures that $\rho$ is zero at
    the origin despite being valid coordinates in a polytope.

  \item The smooth clustering network, as opposed to a standard clustering, also
    ensures differentiability of the network as it will become important for the
    series expansion.

  \item The design of the clustering network ensures that $c(0) = e_{1}$ so
    that values of $\rho$ close to zero are associated with the
    first cluster and that $c(\mu(0)) \otimes \mu(0) = (1, 0) \in
    \R^{\nclstr r}$ represents the target vector $\xs = 0$ exactly.

\end{itemize}

With the following lemma, we summarize the resulting properties when employing
the autoencoder in series expansions of the coefficient function.

\begin{lemma} \label{lem:taylor-of-crhorho}
  Given the smooth clustering network $c$ as defined
  in~\cref{eq:c-clustering-net} and the corresponding LPV coefficient function
  \begin{equation} \label{eq:rho-to-LPV-coeff}
    \rho \to A(\rho) = \tA\Bigr( W\, \bigl(c(\rho) \otimes \rho) \Bigr)
  \end{equation}
  defined through the linear map $\tA\colon \R^{n}\to \R^{n\times n}$ and the
  decoder~\cref{eq:reconstruction-w-W}.
  Then, it holds that
  \begin{itemize}
    \item[(a)] the coefficient function~\cref{eq:rho-to-LPV-coeff} is
      differentiable,
    \item[(b)] the Jacobian of the clustering network $c$ is zero at $\rho = 0$,
      i.e.,
      \begin{equation*}
        c'(\rho)\bigl|_{\rho = 0} = 0 \in \R^{q,r},
      \end{equation*}
      and
    \item[(c)] the second-order Taylor expansion of~\cref{eq:rho-to-LPV-coeff}
      coincides with the first-order expansion and reads
      \begin{equation}\label{eq:Arho-first-second-order}
        A(\rho) \approx A(0) + \sum_{j}^{r} \rho_{j} \tA(w_{j}),
      \end{equation}
      where $w_{j}$ are the first $r$ columns of $W \in \R^{n\times \nclstr r}$
      as in~\cref{eq:reconstruction-w-W}.
  \end{itemize}
\end{lemma}
\begin{proof}
  As for the result (a), we note that, by construction, $c$ is a concatenation
  of arbitrarily often differentiable functions.
  As for the part~(b), by direct computations, we can confirm that
  $c'(\Fclstr(\rho) + \hat{e}_{1}) \big\rvert_{\rho = 0} =
  c'(\hat{e}_{1})= 0 \in \R^{\nclstr \times r}$.

  Since the coefficient map
  \begin{equation*}
    \rho \to \tA\Bigr(W\, \bigl(c(\rho) \otimes \rho)\Bigr)
  \end{equation*}
  is linear in components that concern $W$ and $\tA$ (cf.
  \Cref{eq:reconstruction-w-W} and \Cref{rem:Alinear}), we can reduce the
  argument to the nonlinear part $\rho \to c(\rho) \otimes \rho$.
  With $\rho = 0$, from
  \begin{equation*}
    c(\rho + h) \otimes (\rho + h) - c(\rho) \otimes \rho = c(0 + h) \otimes h,
  \end{equation*}
  we derive the first and second-order expansion terms of the map
  $\rho \to c(\rho) \otimes \rho$ in the increment $h$ at $\rho = 0$ as
  \begin{equation*}
    c(0 + h) \otimes h = \left( c(0) + c'(0) h + \mathcal R(h) \right)
    \otimes h
  \end{equation*}
  with the remainder term $\mathcal R(h)\in o(|h|)$, meaning that
  $\mathcal{R}$ approaches $0$ faster than $h$ itself.
  Thus, the second-order expansion of $\rho \to c(\rho)\otimes \rho$ around the
  origin $0$ reads
  \begin{equation*}
    c(\rho) \otimes \rho \approx c(0) \otimes \rho + (c'(0)\rho) \otimes \rho
      = c(0)\otimes \rho + (c'(0) \otimes I)(\rho \otimes \rho),
  \end{equation*}
  where $c'(0)\in \R^{\nclstr \times r}$ is the Jacobian of $c$ at
  $\rho = 0$.
  With $c'(0) = 0$, the second-order terms vanish so that with
  $c(0) = \hat{e}_{1}$, the second-order expansions reduces to the first order
  approximation of the polytopic reconstruction that reads
  \begin{equation*}
    \tx(\rho) \approx W \left( c(0) \otimes \rho \right) =
      W \left( \hat{e}_{1} \otimes \rho \right) =
      W \begin{bmatrix} \rho \\ 0 \\ \vdots \\ 0 \end{bmatrix} =
      \sum_{j=1}^{r} \rho_{j}  w_{j}.
  \end{equation*}
  Finally, formula~\cref{eq:Arho-first-second-order} follows from the
  linearity of $\tA$ and $W$.
\end{proof}

\begin{remark}
  We note the following consistency property in the smooth clustering.
  By construction, the first $r$ columns of $W$ describe the local basis for
  the reconstruction in the cluster associated with $\rho = 0$.
  Accordingly, the linearization around the origin only considers the basis
  that is associated with the origin.
\end{remark}


\section{Higher order series expansions of LPV approximations in SDRE feedback
  design}%
\label{sec:hoSDRE-expansion}

Although, as shown in \Cref{lem:taylor-of-crhorho}, the second-order
expansion of the coefficient function $A$ coincides with the expansion of first
order, the feedback law will still exhibit second-order dependencies as we will
lay out in this section.
For the \emph{quasi}-LPV system~\cref{eq:nonlsys-lpv}, the nonlinear SDRE feedback
design defines the input $u$ as
\begin{equation} \label{eq:sdre_fblaw}
  u(t) = -B^{\trans} P(\rho(x(t))) x(t),
\end{equation}
where $P(\rho)$ solves the SDRE
\begin{equation} \label{eq:sdre}
  A(\rho)^{\trans} P(\rho) + P(\rho) A(\rho) - P(\rho) B B^{\trans} P(\rho)
    = -C^{\trans} C.
\end{equation}
We note that for \emph{quasi}-LPV systems the parameter changes continuously
with the states.
Since the repeated solve of~\cref{eq:sdre} is costly, series expansions of the
feedback law~\cref{eq:sdre_fblaw} have been proposed.
We now briefly derive the formulas for higher order expansions of the solution
to~\cref{eq:sdre} with nonlinear parameter dependencies in $A$.
For linear parameter dependencies, the first-order approximation has been
developed in~\cite{BeeTB00} for a single parameter dependency
and has been extended in~\cite{AllKS23} to the multivariate case.

If the state dependency of the SDRE solution is parametrized through $\rho$,
the multivariate Taylor expansion of $P$ about $\rho_{0} = 0$ up to order
$\xrdr$ reads
\begin{equation} \label{eqn:Pexp}
  P(\rho) \approx P_{0} + \sum_{1 \leq \lvert \alpha \rvert \leq \xrdr}
    \rho^{(\alpha)} P_{\alpha},
\end{equation}
where $\alpha =(\alpha_{1}, \ldots, \alpha_{r}) \in \N^{r}$ is a multiindex with
$\lvert \alpha \rvert := \sum_{i=1}^{r} \alpha_{i}$, where the coefficients are
$\rho^{(\alpha)} := \rho_{1}^{\alpha_{1}} \rho_{2}^{\alpha_{2}} \cdots
\rho_{r}^{\alpha_{r}}$, and where, importantly, all $P_{\alpha}$ are constant
symmetric matrices, namely
\begin{equation*}
  P_{\alpha} =  \frac{1}{\alpha_{1}! \alpha_{2}! \cdots \alpha_{r}!}
    \frac{\partial ^{\lvert \alpha \rvert}}{\partial_{\rho_{1}}^{\alpha_{1}}
    \partial_{\rho_{2}}^{\alpha_{2}} \cdots \partial_{\rho_{r}}^{\alpha_{r}}}
    P(0).
\end{equation*}
We note that, because the function $P$ is not known, this formula is merely stated
to confirm that the $P_\alpha$ are constant matrices. Practically, in what follows, we
will compute the $P_\alpha$ directly by equating coefficients in truncated
series.

For the nonlinear parameter dependencies in $A(\rho)$, we will consider a
truncated series expansion of $A$:
\begin{equation} \label{eqn:Aexp}
  A(\rho) \approx A_{0} + \sum_{1\leq \lvert \alpha \rvert \leq \xrdr}
    \rho^{(\alpha)}A_{\alpha},
\end{equation}
where $A_{\alpha}$ are constant matrices; see~\cref{eq:Arho-first-second-order}
that describes the first and second-order expansions as follows:
For $p = 1$, we can identify all
$\lvert \alpha \rvert=p$ with $k=1,\dotsc, r$ and obtain
$\rho^{(\alpha)}A_{\alpha}=\rho_{k} \widetilde{A}(w_{k})$,
where $w_{k}$ is the $k$-th column of $W$.
Following \Cref{lem:taylor-of-crhorho}, in the considered case with
the structure of the nonlinear decoding as $\rho \to W (c(\rho) \otimes \rho)$,
the second-order expansion coincides with the first-order one.
We also note that, if $\widetilde{A}$ is known (see \Cref{rem:Alinear}),
then $A_{0}$ is readily inferred without knowledge of $\rho$ or $W$.

We insert the expansions of $A$ and $P$, namely~\cref{eqn:Pexp,eqn:Aexp}, into
the SDRE~\cref{eq:sdre} and solve for the coefficients $P_{\alpha}$ by matching
the individual powers of $\rho$.
For $\alpha = 0$, we obtain the classical standard Riccati equation
\begin{equation} \label{eqn:riccati}
  A_{0}^{\trans} P_{0} + P_{0} A_{0} - P_{0} B B^{\trans} P_{0} = -C^{\trans} C,
\end{equation}
as it is used for the design of feedback controllers for linear systems.
Moving on, for $\lvert \alpha \rvert = 1$, the coefficients $P_{\alpha}$ are
obtained from requiring that the first-order terms fulfill
\begin{equation*}
  \begin{split}
    P_{0} \left(\sum_{\lvert \alpha \rvert =1} \rho^{(\alpha)}
      A_{\alpha} \right) + \left(\sum_{\lvert \alpha \rvert =1} \rho^{(\alpha)}
      A_{\alpha} \right)^{\trans} P_{0} -
      \sum_{\lvert \alpha \rvert =1} \rho^{(\alpha)} P_{\alpha} B
      B^{\trans} P_{0} - P_{0} B B^{\trans}
      \sum_{\lvert \alpha \rvert =1} \rho^{(\alpha)} P_{\alpha} = 0,
  \end{split}
\end{equation*}
for arbitrary values of $\rho$.
This defines the $P_{\alpha}$ as the solutions to the Lyapunov equations
\begin{equation} \label{eqn:lyap1}
  \left( A_{0} - B B^{\trans} P_{0} \right)^{\trans} P_{\alpha} +
    P_{\alpha} \left( A_{0} - B B^{\trans} P_{0} \right) =
    -A_{\alpha} ^{\trans} P_{0} - P_{0} A_{\alpha},
\end{equation}
for the $r$ multiindices in $\alpha \in \N^{r}$ for which
$\lvert \alpha \rvert = 1$ holds; see also~\cite{AllKS23, morHeiW23} where the
same formulas have been derived.
Similarly, we may obtain the $P_{\alpha}$ by setting the terms in
which $\rho$ appears quadratic to zero:
\begin{equation*}
  \begin{aligned}
    0 & = P_{0} \left(\sum_{\lvert \alpha \rvert = 2} \rho^{(\alpha)}
      A_{\alpha} \right) + \left( \sum_{\lvert \alpha \rvert = 1}
      \rho^{(\alpha)} P_{\alpha} \right) \left(\sum_{\lvert \alpha \rvert = 1}
      \rho^{(\alpha)} A_{\alpha} \right) +
      \left(\sum_{\lvert \alpha \rvert = 2} \rho^{(\alpha)} P_{\alpha}
      \right) A_{0} + \dotsc \\
    & \qquad {}-{}
      \left(\sum_{\lvert \alpha \rvert =1} \rho^{(\alpha)} P_{\alpha} \right)
      B B^{\trans} \left(\sum_{\lvert \alpha \rvert = 1} \rho^{(\alpha)}
      P_{\alpha} \right) -
      \left( \sum_{\lvert \alpha \rvert = 2} \rho^{(\alpha)} P_{\alpha}
      \right) B B^{\trans} P_{0} - \dotsc,
  \end{aligned}
\end{equation*}
where the terms omitted via $\dotsc$ stand for the transposed counterparts.
Let $\astr$ be a multiindex with $\lvert \astr \rvert = 2$ and let
$\astro$ and $\astrt$ be those two multiindices with $\lvert \astro \rvert =
\lvert \astrt \rvert = 1$ that satisfy $\astro + \astrt = \astr$.
Then, the coefficient $P_{\astr}$ that cancels the corresponding contributions
in the equation above is defined as the solution to the Lyapunov equation
\begin{equation} \label{eqn:lyap2}
  \begin{aligned}
    & (A_{0} - B B^{\trans} P_{0})^{\trans} P_{\astr} +
      P_{\astr} (A_{0} - B B^{\trans} P_{0})\\
    & = - P_{0} A_{\astr} - A_{\astr}^{\trans} P_{0} - P_{\astro} A_{\astrt} -
      A_{\astrt}^{\trans} P_{\astro} - P_{\astrt} A_{\astro}\\
    & \quad{}-{}
      A_{\astro}^{\trans} P_{\astrt} + P_{\astro} BB^{\trans} P_{\astrt} +
      P_{\astrt} BB^{\trans} P_{\astro}.
  \end{aligned}
\end{equation}
Note that for given $r$ and $\xrdr$, there exist exactly
$\binom{r+\xrdr-1}{\xrdr}$ unique such multiindices $\alpha$ such that
$\rvert \alpha \lvert = \xrdr$.

Finally, we note that for such a series expansion of $P$, the feedback
law~\cref{eq:sdre_fblaw} reads
\begin{equation}\label{eq:sdre_fblaw_xpanded}
  u(t)= -B^{\trans} \left( \sum_{\rvert \alpha \lvert \leq p} \rho^{(\alpha)}
    P_{\alpha} \right) x(t) =:  u_{P_{\lvert \alpha \rvert \leq p}}(t),
\end{equation}
which with the precomputed $B^{\trans} P_{\alpha}$ is efficiently evaluated by
adapting the values of $\rho$.
We also note that for $p = 0$, this formula reduces to the standard \emph{linear
quadratic regulator} (LQR) feedback with respect to a linearization about
$\xs = 0$.


\section{Solving the high-dimensional matrix equations}%
\label{sec:highdimme}

To compute the solutions to the occurring high-dimensional matrix
equations~\cref{eqn:riccati,eqn:lyap1,eqn:lyap2},
some reformulations and state-of-the-art approaches are needed.
In the considered setup, we can exploit that the ranks of $B$ and $C$ are much
smaller than the state-space dimension $n$ and that $A_{0}$ in~\cref{eqn:Aexp}
is sparsely populated, so that
the first equation~\cref{eqn:riccati}---a standard Riccati
equation---can be efficiently solved for a low-rank approximation of the
symmetric positive semi-definite stabilizing solution
$Z_{0} Z_{0}^{\trans} \approx P_{0}$, with
$Z_{0} \in \R^{n \times k}$ and $k \ll n$; see~\cite{BenB16, Sti18}, e.g.,
via the low-rank Newton-Kleinman iteration~\cite{BenHSetal16, SaaW24},
the incremental subspace iteration~\cite{LinS15}
or the RADI method~\cite{BenBKetal18, BerF24}.

However, for the two types of Lyapunov equations~\cref{eqn:lyap1}
and~\cref{eqn:lyap2} occurring through the series expansion, some
reformulations are needed to compute the solutions in this large-scale setting;
namely, the right-hand sides need to be reformulated into symmetric low-rank
factorizations.
Given the approximation $Z_{0} Z_{0}^{\trans} \approx P_{0}$ of the initial
Riccati equation~\cref{eqn:riccati}, the right-hand side in~\cref{eqn:lyap1}
can be rewritten as
\begin{equation} \label{eqn:lyap1rhs}
  -A_{\alpha}^{\trans} P_{0} - P_{0} A_{\alpha} \approx
    -\begin{bmatrix} A_{\alpha}^{\trans} Z_{0} & Z_{0} \end{bmatrix}
    J
    \begin{bmatrix} Z_{0}^{\trans} A_{\alpha} \\ Z_{0}^{\trans} \end{bmatrix},
\end{equation}
where the symmetric yet indefinite center matrix is given by
$J = \begin{bmatrix} 0 & I_{k} \\ I_{k} & 0 \end{bmatrix}$, with $I_{k}$
denoting the $k$-dimensional identity matrix and $k$ being the number of columns
of $Z_{0}$.
For Lyapunov equations~\cref{eqn:lyap1} with such symmetric factorized
right-hand sides~\cref{eqn:lyap1rhs}, the ADI approach has been extended to
compute appropriate solution approximations of the form
$L_{\alpha} D_{\alpha} L_{\alpha}^{\trans} \approx P_{\alpha}$ following the
format of the right-hand side; see~\cite{LanMS14}.
Furthermore, while the coefficient matrix in~\cref{eqn:lyap1}, namely
\begin{equation} \label{eqn:Alowrank}
  A - B B^{\trans} P_{0} = A - B V^{\trans}
\end{equation}
with $V = B^{\trans} P_{0} \in \R^{n \times m}$, is dense and high-dimensional,
it is not necessary to assemble this matrix at any point
because~\cref{eqn:Alowrank} is a large-scale sparse matrix with low-rank
perturbation.
Solution methods for large-scale Lyapunov equations only rely on
matrix-vector operations or the solution to linear systems
with~\cref{eqn:Alowrank}.
The latter can be efficiently computed using the Sherman-Morrison-Woodbury
matrix inversion formula or augmented matrix approach~\cite{GolV13} such that
only the solution of linear systems with large-scale sparse coefficient
matrices and matrix-vector products are needed.
It is important to observe that the coefficient~\cref{eqn:Alowrank} is Hurwitz,
i.e., all eigenvalues lie in the left open half-plane, such that the established
solution methods can be applied.

Similarly, we can derive for the right-hand side of~\cref{eqn:lyap2}
symmetric indefinite low-rank factorizations of the form $L D L^{\trans}$.
Thereby, the factor are
\begin{equation} \label{eqn:lyap2rhs1}
  L^{\trans} = \begin{bmatrix} Z_{0}^{\trans} A_{\astr} \\
    Z_{0}^{\trans} \\
    L_{\astro}^{\trans} (A_{\astrt} - B K_{\astrt}) \\
    L_{\astrt}^{\trans} (A_{\astro} - B K_{\astro})\\
    L_{\astro}^{\trans} \\
    L_{\astrt}^{\trans} \end{bmatrix},
\end{equation}
and
\begin{equation} \label{eqn:lyap2rhs2}
  D = \begin{bmatrix} J & 0 & 0 & 0 & 0 \\ 0 & 0 & 0 & D_{\astro} & 0 \\
    0 & 0 & 0 & 0 & D_{\astrt} \\
    0 & D_{\astro} & 0 & 0 & 0 \\ 0 & 0 & D_{\astrt} & 0 & 0 \end{bmatrix},
\end{equation}
where we have the approximations of the solutions to the Lyapunov
equations~\cref{eqn:lyap1} as
$L_{\astro} D_{\astro} L_{\astro}^{\trans} \approx P_{\astro}$ and
$L_{\astrt} D_{\astrt} L_{\astrt}^{\trans} \approx P_{\astrt}$,
and the intermediate feedback matrices $K_{\astro} = B^{\trans} P_{\astro}$ and
$K_{\astrt} = B^{\trans} P_{\astrt}$.
As for~\cref{eqn:lyap1}, this reformulation of the constant terms
of~\cref{eqn:lyap2} allows the application of the same computational
procedures for the solution, for example, the $LDL^{\trans}$-factorized ADI
method~\cite{LanMS14}.
The coefficient matrix in~\cref{eqn:lyap2} is the same as in~\cref{eqn:lyap1}
such that it has the same properties being sparse plus low-rank and Hurwitz.
Therefore, we can apply the same numerical approaches to solve~\cref{eqn:lyap2}
as for~\cref{eqn:lyap1}.

In both cases of Lyapunov equations~\cref{eqn:lyap1,eqn:lyap2}, we have that
the inner dimensions of the reformulated right-hand
side terms~\cref{eqn:lyap1rhs,eqn:lyap2rhs1,eqn:lyap2rhs2} depend on dimensions
of the factors of the solution approximations to~\cref{eqn:riccati,eqn:lyap1}.
Therefore, we can expect the dimensions of the right-hand side factors to be
larger than the dimensions of the input and output matrices used in the Riccati
equation~\cref{eqn:riccati}.
In general, we strongly recommend to compress the right-hand side factors at
least to their numerical rank using economy-size QR and eigenvalue
decompositions.
All right-hand sides of Lyapunov equations considered in this paper can be
written as indefinite symmetric factorizations of the form $LDL^{\trans}$, with
$L \in \R^{n \times k}$ and $D \in \R^{k \times k}$ symmetric.
For the compression, first, compute the economy-size QR of $L = QR$, with
$Q \in \R^{n \times k}$ and $R \in \R^{k \times k}$.
Next, compute the eigenvalue decomposition of the symmetric matrix
\begin{equation} \label{eqn:symprod}
  R D R^{\trans} = U \Lambda U^{\trans},
\end{equation}
where $\Lambda \in \R^{k \times k}$ is a diagonal matrix of real eigenvalues
and $U \in \R^{k \times k}$ is the matrix of the corresponding orthogonal
eigenvectors.
Next, the numerical rank is determined based on the diagonal matrix $\Lambda$.
A typical choice here is to truncate all eigenvalues $\lambda_{i}$ that are
small in absolute value, e.g.,
\begin{equation*}
  \lvert \lambda_{i} \rvert < k \cdot \max_{i = 1, \dots, k}
    \lvert \lambda_{i} \rvert \cdot \epsilon,
\end{equation*}
where $\epsilon$ is the machine precision.
Note that the original matrix $D$ is assumed to be indefinite such that large
positive as well as negative eigenvalues may occur, which need to be kept in
the compression.
For simplicity, assume that the eigenvalues in $\Lambda$ are ordered by
absolute size such that we have
\begin{equation*}
  \Lambda = \begin{bmatrix} \Lambda_{1} & 0 \\ 0 & \Lambda_{2} \end{bmatrix}
  \quad\text{and}\quad
  U = \begin{bmatrix} U_{1} & U_{2} \end{bmatrix},
\end{equation*}
where $\Lambda_{1}, U_{1}$ contain the eigenvalues and eigenvectors that remain
and $\Lambda_{2}, U_{2}$ contain the eigenvalues and eigenvectors that are
truncated.
This results in the compression
\begin{equation*}
  L D L^{\trans} \approx Q U_{1} \Lambda_{1} U_{1}^{\trans} Q^{\trans} =
    \widetilde{L} \widetilde{D} \widetilde{L}^{\trans},
\end{equation*}
where $\widetilde{L} = Q U_{1}$ and $\widetilde{D} = \Lambda_{1}$, which is
used to compute the solution to the matrix equations~\cref{eqn:lyap1,eqn:lyap2}
instead of the original right-hand sides.


\section{Numerical experiments}%
\label{sec:num-exps}

The code, data and results of the numerical experiments presented here can be
found at~\cite{supHeiKW24}.


\subsection{Numerical setup}

\begin{figure}[t]
  \centering
  \begin{minipage}{.49\linewidth}
    \centering
    \includegraphics[width = \linewidth]{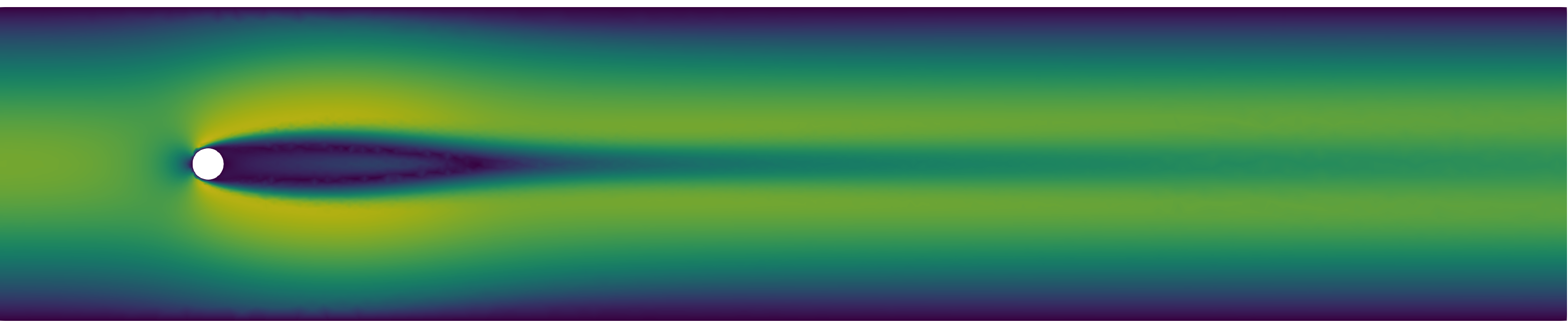}
  \end{minipage}%
  \hfill%
  \begin{minipage}{.49\linewidth}
    \centering
    \includegraphics[width = \linewidth]%
      {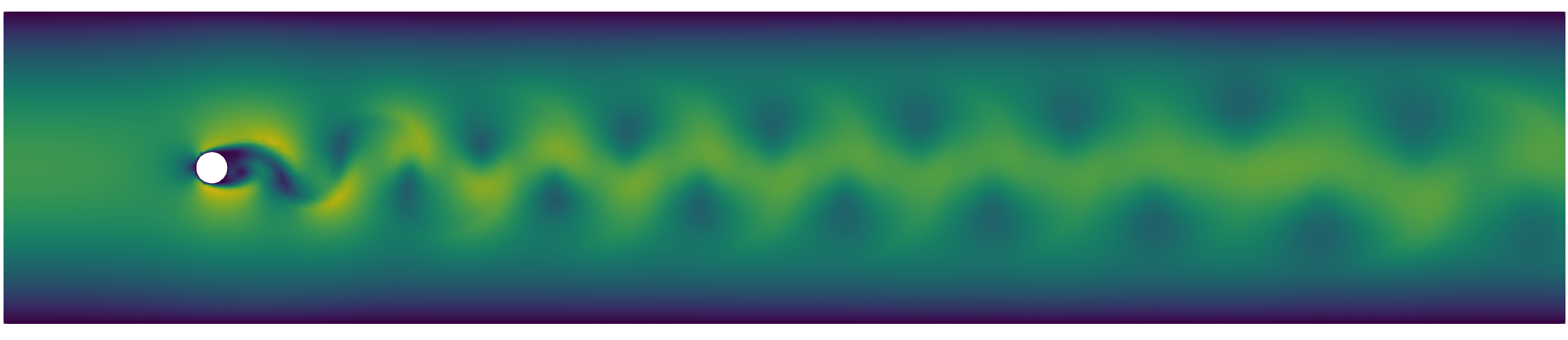}
  \end{minipage}

  \caption{Spatial domain and snapshots of the flow behind a cylinder in the
    unstable steady state (left) and in the periodic vortex shedding regime
    (right); cf.~\cite{morHeiW23}.}
  \label{fig:setup}
\end{figure}

For comparability reasons, we consider here the same problem, numerical setup,
and the same data points for training the models as in~\cite{morHeiW23}.
In short, the flow control problem of stabilizing the unstable steady-state of
the flow behind a cylinder is modeled through the difference of the velocity $v$
to the steady-stated and observed through outputs $y = C x$.
It is determined by the semi-discrete Navier-Stokes equations
\begin{equation} \label{eq:nse-projected}
  M \dot x = N(x) x  + Bu, \quad y=Cx, \quad x(0)= 0 \in \R^{n},
\end{equation}
written in so-called \emph{divergence-free} coordinates of the difference state
$x$, with $N\colon x(t)\to N(x(t))\in \R^{n\times n}$ modeling the discrete
diffusion-convection operator and $C \in \R^{6\times n}$ and
$B \in \R^{n\times 2}$ the discrete observation and control operators;
cf.~\cite{morHeiW23}.
See also \Cref{fig:setup} for a picture of the spatial domain and two snapshots
of the flow behind a cylinder showing the flow velocities in the unstable
steady state and with the periodic vortex shedding behavior.
In the considered setup, the \emph{Reynolds number} is set to~$60$ and the
dimension of the spatially discretized model to~$n = 51\,194$.
For the time integration, a semi-explicit Euler scheme is employed, which
treats the feedback and the nonlinearity explicitly.

\begin{table}[t]
  \centering
  \caption{Model complexity of the different PAE schemes and POD:
    \# stands for \emph{number of}.
    We distinguish between linear (lin.) and nonlinear (nonl.)
    mappings in the encoding and decoding.}
  \label{tab:pae}
  \vspace{.5\baselineskip}

  \begin{tabular}{c|c|c|c|c} \hline
    Hyperparameters \textbackslash~ \texttt{Scheme} & \POD 5 & \PAEqr 15 &
    \PAEqr 35 & \makecell{ \PAEqr 35 \\ \texttt{(p=1)}}\\
    \hline
    reduced dimension $r$ &5 & 5 & 5 & 5\\
    \hline
    \#clusters $\nclstr$&1  & 1 & 3 & 3\\
    \hline
    \#encoding parameters& 255970 & 42810 & 42810 & 42810 \\
    \hline
    \#encoding layers &1(lin.) & 16(nonl.) & 16(nonl.) & 16(nonl.)\\
    \hline
    \#decoding parameters& 255970 & 255970 & 15+767910 & 255970 \\
    \hline
    \#decoding layers &1& 1 & 1(nonl.)+1(lin.) & 1\\
    \hline
  \end{tabular}
\end{table}

We consider the steady-state solution as target state and aim for the
stabilization of the difference system~\cref{eq:nse-projected}.
To trigger instabilities and to collect training data, a
test input
\begin{equation} \label{eq:testu}
  u(t) = \begin{bmatrix} \sin(t) & 0 \end{bmatrix}^{\trans},
\end{equation}
is applied.
We choose a deep convolutional encoder with depthwise separable convolutions to
cope with the high computational complexity and memory requirements in this
rather high-dimensional setup; see~\cite{HeiK24} for details.
As a result the model complexity in terms of parameters that define the encoder
and decoder can be lower than for, say, a comparable POD reduction;
see \Cref{tab:pae}.
For training the PAEs, we adopt a semi-supervised learning methodology
with the loss function
\begin{equation*}
  \mathcal{L}(v^{(i)}) := \lambda \| \tx^{(i)} - x^{(i)}\|_{M} -
    \ell(\rho^{(i)}) \cdot \log(c(\rho^{(i)})),
\end{equation*}
where $x^{(i)}$ is the $i$-th data point,
$\rho^{(i)} := \mu(x^{(i)})$ and
$\tx^{(i)} :=\nu(\rho^{(i)})$
for the corresponding encoder and decoder $\mu$ and $\nu$,
where $\ell$ assigns labels that were precomputed by $k$-means clustering in
$\R^{r}$, and where $\lambda > 0$ is a weighting factor between the
reconstruction error (measured in the $M$-weighted $L^{2}$-norm) and
the clustering error.
We note that the chosen \emph{cross-entropy} loss is commonly used to
smoothly assign labels and that the labeling $\ell$ is defined such that the
target $\rho^{\ast} = 0$ is assigned to the first cluster.
As it is standard practice, the loss is computed as the average over a
\emph{batch} of data points.
Here, we consider $401$ data points equally distributed in the
time interval $[0, 0.5]$ as training data, and we use the
\emph{Adam optimizer} with learning rate $0.005$, $3200$ epochs, a
batch size of $64$ and the loss weight $\lambda = 100$.


\subsection{Reconstruction errors}

\begin{table}[t]
  \centering
  \caption{Average of the reconstruction errors (scaled with $10^{3}$) of \POD{r} and
    \PAEqr{q}{r} on the time range
    $[0,0.5]$ used for training.}
  \label{tab:avrgd_rec_errors}
  \vspace{.5\baselineskip}

  \begin{tabular}{c|c|c|c|c|c|c|c} \hline
    \texttt{Scheme \textbackslash ~ r} & 2 & 3 & 4 & 5 & 6 & 7 & 8\\
    \hline
    \POD r & 1.340& 0.610& 0.305 & 0.149 & 0.076 & 0.038 & 0.020\\
    \hline
    \PAEqr 1r  & 1.335 & 0.614 & 0.326 & 0.181 & 0.118 & 0.121 & 0.087\\
    \hline
    \PAEqr 3r  & 0.189 & 0.135 & 0.104 & 0.066 & 0.075 & 0.097 & 0.079 \\
    \hline
    \PAEqrp 3r1 & 17.670 & 7.501 & 5.217 & 5.628 & 2.745 & 6.794 & 2.307\\
    \hline
  \end{tabular}
\end{table}

\begin{figure}[t]
  \centering
  \resizebox{.99\linewidth}{!}{\input{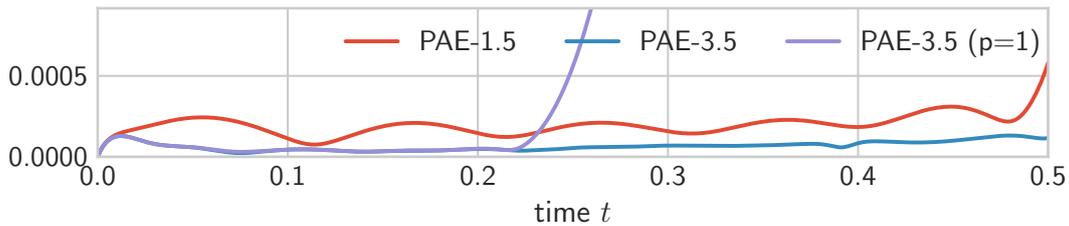}}
  \vspace{-\baselineskip}

  \caption{Trajectory with pointwise reconstruction errors
    $\|x(t^{(i)}) - \tx(t^{(i)})\|_{M}$.
    The improvement of using $3$ clusters is clearly visible as is the initial
    alignment with the nonlinear approach \PAEqr{3}{5}
    (and the subsequent deviation from it) of the first-order expansion.}
  \label{fig:rec-err}
\end{figure}

We evaluate the reconstruction performance of the PAE in the averaged $M$-norm
error
\begin{equation*}
  \frac{1}{N}\sum_{i=1}^{N} \|x^{(i)} - \tx^{(i)} \|_{M},
\end{equation*}
where $x^{(i)}$ is the $i$-th data point and $\tx^{(i)}=\nu(\mu(x^{(i)}))$
for the corresponding encoder and decoder $\mu$ and $\nu$.
As laid out in \Cref{tab:avrgd_rec_errors}, the \PAEqr{1}{r} is comparable
to the standard linear approach of POD for low dimensions $r = 2, 3, 4$, whereas
the nonlinear encoding through the cluster selection significantly
improves the reconstruction in low dimensions.
As for the series expansion of the control law, only the first-order
approximation is relevant, we tabulate this reconstruction error as well.
Note that the bad average performance is due to the strong deviation for larger
$t$, i.e., larger values of $\rho$, whereas for small $t$, the first-order
approximation closely aligns with the nonlinear approach; see
\Cref{fig:rec-err}.

\begin{figure}[t]
  \centering
  \resizebox{.99\linewidth}{!}{\input{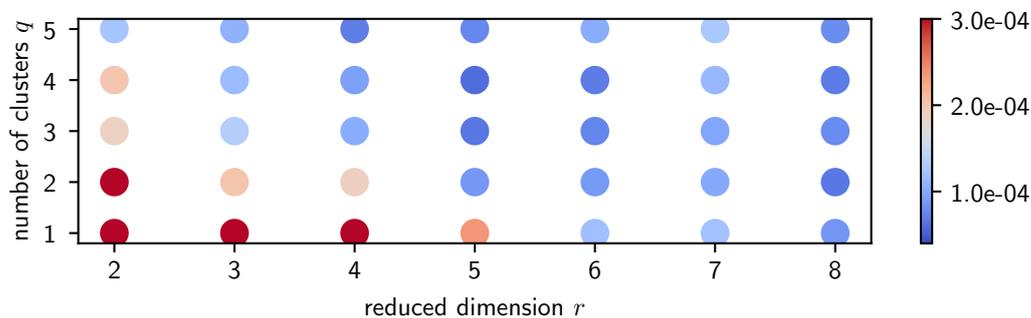}}
  \vspace{-\baselineskip}

  \caption{Grid search for the selection of $r$ and $\nclstr$ with lowest
    averaged reconstruction error.
    A trade-off between the number of clusters and the reduced dimension
    becomes clearly visible.}
  \label{fig:grid-search}
\end{figure}

Before turning to the controller design, we use a grid search to identify a
suitable parameter setup for the reduced-order LPV approximation.
As shown in \Cref{fig:grid-search}, the pair $(q, r) = (4, 5)$ achieves the
lowest reconstruction error in the considered domain.
Nonetheless, we selected $(q, r) =(3,5)$ as a compromise of both accuracy and
complexity.


\subsection{SDRE approximation through series expansions}

\begin{figure*}[t]
  \centering
  \resizebox{.99\linewidth}{!}{\input{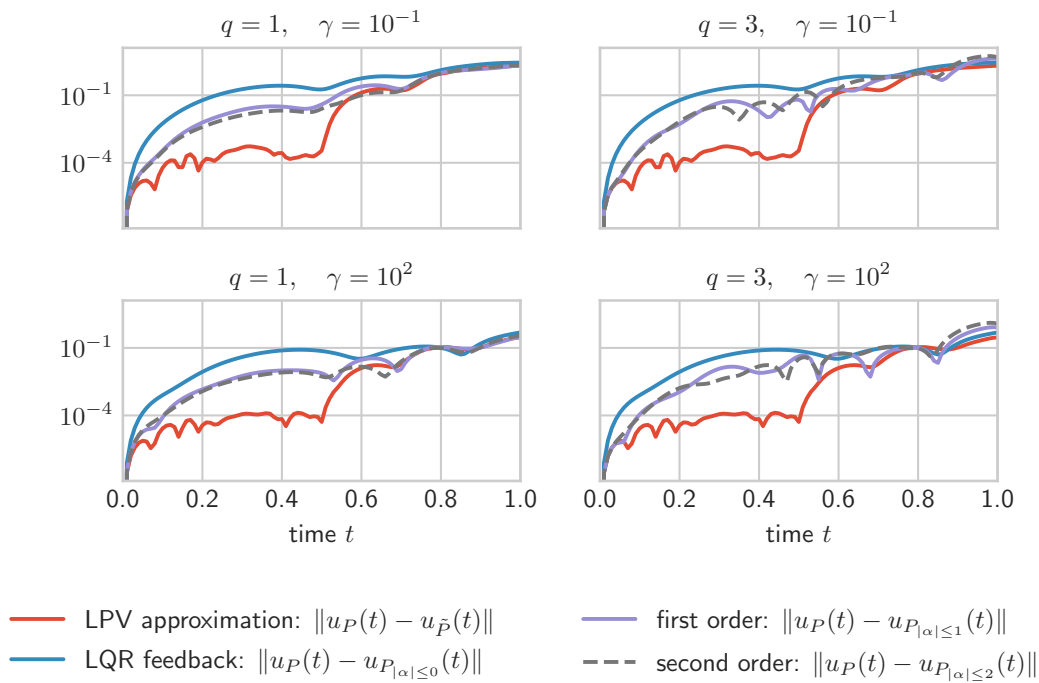}}
  \vspace{-.5\baselineskip}

  \caption{Differences to the exact SDRE feedback along a sample trajectory for
    feedback approximations through the LPV approximation of the nonlinearity by
    \PAEqr{q}{5}, for $q = 1, 3$, and the corresponding zeroth, first and
    second-order series expansions of the SDRE feedback.
    The LPV approximation provides the most accurate approximation to the true
    SDRE feedback in its training regime up to $0.5$.
    The first and second-order approximations also strongly improve over the
    classical LQR feedback.}
  \label{fig:sdre-approx}
\end{figure*}

We check the approximation of the SDRE feedback along the example trajectory $x$
generated by the test input~\cref{eq:testu}.
For that we compute the true SDRE feedback
\begin{equation*}
  u_{P}(t) = -B^{\trans} P(x(t)) x(t),
\end{equation*}
where $P$ solves the SDRE, and compare
to $u_{\tilde{P}}$ and $u_{P_{|\alpha| \leq p}}$ that denote the feedback
computed through the LPV approximation in~\cref{eq:sdre}
(via \PAEqr{q}{5}, $q = 1, 3$) and the expansions~\cref{eq:sdre_fblaw_xpanded}
for $p = 0, 1, 2$, respectively.
The data in \Cref{fig:sdre-approx} shows that the LPV approximation works
particularly well for $t \leq 0.5$, which is the training regime, and that the
first and second-order expansions significantly improve the zero-order
approximation, which would be the LQR gain, but with a slight effect of the
additional second-order terms.


\subsection{Feedback performance}

\begin{figure}[t]
  \centering
  \resizebox{.99\linewidth}{!}{\input{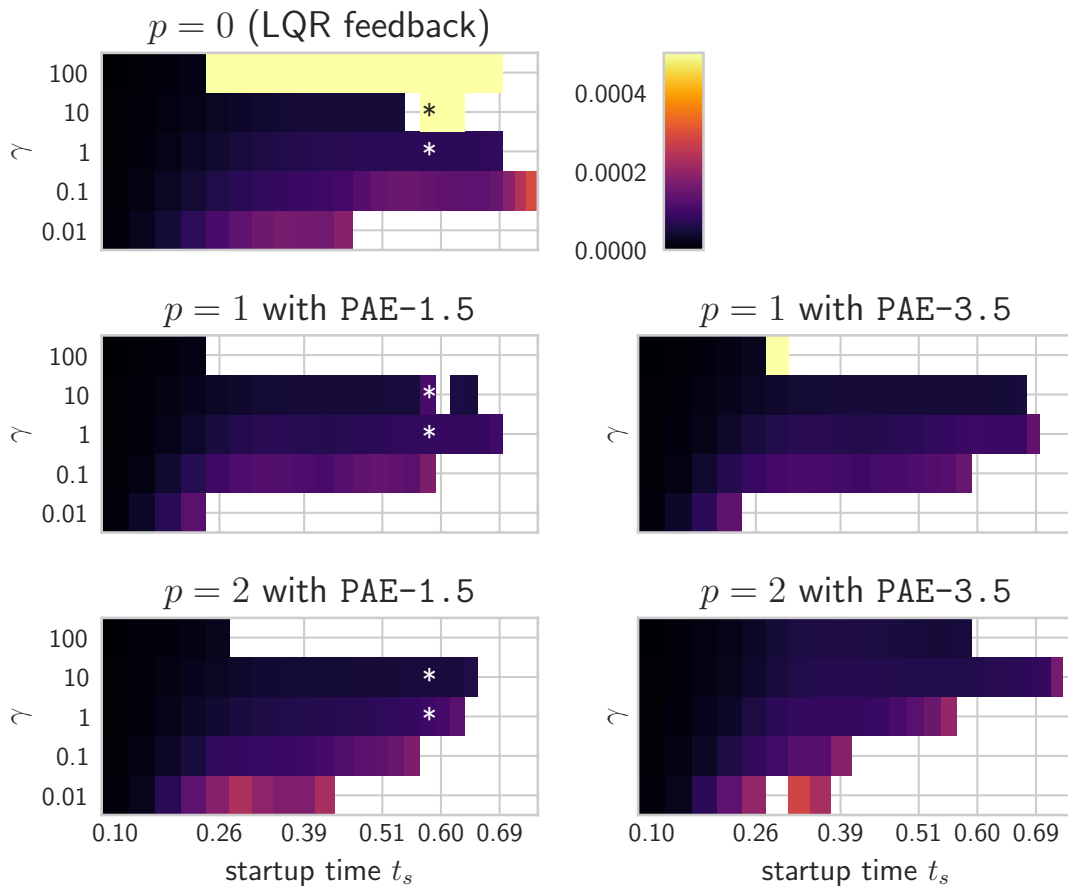}}
  \vspace{-.5\baselineskip}

  \caption{Performance index~\cref{eq:fb-performance-index} map of the linear
    and nonlinear feedbacks for PAE reduced-order models with $r=5$,
    $q = 1, 3$, the SDRE feedback law expansion~\cref{eq:sdre_fblaw_xpanded}
    for $p = 0, 1, 2$ and varying penalization parameter $\gamma$ and
    startup times $t_s$.
    The darker the color, the smaller (i.e., the better) the performance index.
    At the blank spaces, the system blew up.
    The trajectories of the data marked with * are plotted in
    \Cref{fig:fb-trjcheck}.}
  \label{fig:fb-performance}
\end{figure}

\begin{figure*}[t]
  \centering
  \resizebox{.99\linewidth}{!}{\input{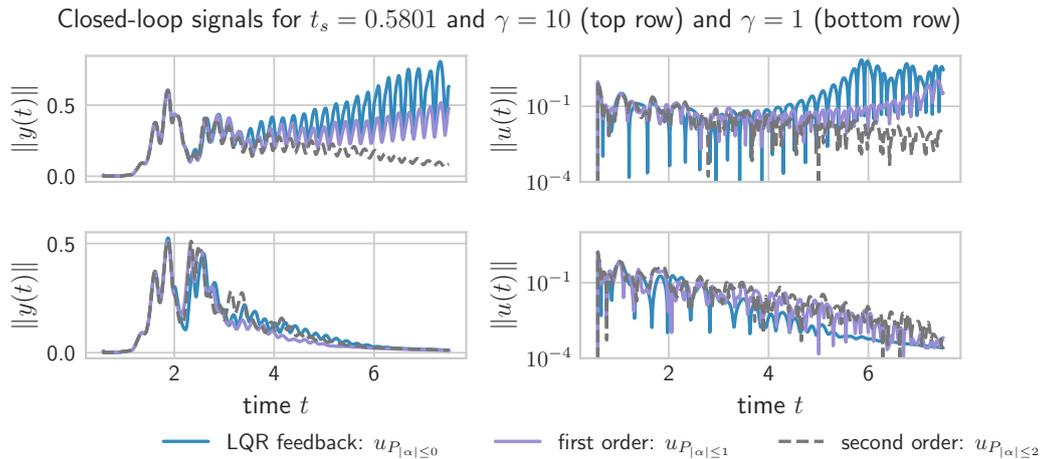}}
  \vspace{-.5\baselineskip}

  \caption{Measured output and feedback input signals for the distinguished data
    points from \Cref{fig:fb-performance}.
    The second-order LPV approximation is capable of stabilizing the nonlinear
    system in situations, in which the first-order approximation begins to
    fail.}
  \label{fig:fb-trjcheck}
\end{figure*}

To evaluate the performance of the approximations for the controller design, we
simulate the closed loop system with varying the two parameters
$\gamma \in \{100, 10, 1, 0.1, 0.001\}$ that stands for the penalization factor
of $u$ in the underlying quadratic cost functional and the startup time
$t_{s} > 0$ that marks the time after which the test input is turned off and
the feedback is applied.
In other words, $\gamma$ is a parameter that weighs on the magnitude of the
control and $t_{s}$ regulates how far the system is deferred from the target
state when the controller is activated.
In \Cref{fig:fb-performance}, for different orders $p = 0, 1, 2$ of
approximation, we plot the averaged measured feedback magnitude
\begin{equation}\label{eq:fb-performance-index}
  \frac{1}{\te} \Big( \int_{t_{s}}^{\te} \| B^{\trans}
    \bigl( \sum_{|\alpha| \leq p} \rho(x(t))^{(\alpha)} P_{\alpha}
    \bigr) x(t) \|^{2} \operatorname{d}t \Big)^{\frac{1}{2}},
\end{equation}
that tends to $0$ for $\te \to \infty$ if stabilization is achieved, and that
is finite if and only if the simulation has not blown up before a finite end
time $\te$.

We plot the values of the performance index~\cref{eq:fb-performance-index} for
$\te = 7.5$ and for
various approximations to the nonlinear SDRE feedback in
\Cref{fig:fb-performance} and the trajectories of selected cases in
\Cref{fig:fb-trjcheck}.
As we have observed before, the $p = 0$ expansion
in~\cref{eq:sdre_fblaw_xpanded}, i.e., the LQR feedback, is remarkably
performant for the optimal choice of the \emph{Tikhonov} parameter $\gamma$.
However, considering higher order expansions with $p = 1, 2$ together
with the option to adapt the underlying approximation scheme, significantly
widens the domain of attraction and can ensure reliable performance in
particular for larger values of $\gamma$; see, e.g., \Cref{fig:fb-trjcheck}
(top row) for the trajectory data for a setup with $\gamma = 10$, where only
the second-order expansion achieves stabilization.
The reason might be that a larger $\gamma$ seeks for small (in magnitude)
control actions, which makes better adaptation pay off while it reduces
overshoots that may destabilize the nonlinear controller.


\section{Conclusion}%
\label{sec:conclusions}

In this work, we considered the use of deep polytopic autoencoders and
state-dependent Riccati equations for the design of nonlinear controllers.
We use the concept of linear parameter-varying approximations of nonlinear
dynamical systems that are parametrized by deep polytopic autoencoders, and we
show that by the particular choice of autoencoders, the parametrization is
differentiable and has vanishing second-order terms.
For the nonlinear feedback design, we extend the series expansion of
state-dependent Riccati equations to second-order terms for which we provide
computable expressions in form of indefinite Lyapunov and Riccati equations.
Finally, we test the proposed framework on the classical control example
modeling the flow past a cylinder.

The proposed nonlinear autoencoder with smooth clustering significantly
improves the parametrization quality at very low dimensions.
For controller design, we established that for second-order expansions of the
state-dependent feedback law, only the first-order expansions of the
reduced-order linear parameter-varying coefficient are relevant, which reduces
the computational effort by large.
Regarding the performance of the resulting feedback, the overall picture is
diffuse.
While the linear quadratic regulator approach, for suitable weight parameters,
ensured stability also well away from the origin, the nonlinear feedback design
based on different approximation schemes, clearly extended the basin of
attraction.
As the choice of the approximation was solely based on a single open-loop
simulation in a time range that did not cover the full range of the
closed-loop simulation, in future work, we expect to gain additional
performance by optimizing the approximation in view of feedback control.


\section*{Acknowledgments}%
\addcontentsline{toc}{section}{Acknowledgments}

Parts of this work have been carried out while Jan Heiland was with the
Max Planck Institute for Complex Technical Systems, Magdeburg, and with the Otto
von Guericke University Magdeburg.

Jan Heiland and Yongho Kim were supported by the German Research Foundation
(DFG) Research Training Group 2297 ``Mathematical Complexity
Reduction (MathCoRe)'', Magdeburg.
Steffen W. R. Werner acknowledges Advanced Research Computing (ARC) at Virginia
Tech for providing computational resources and technical support that have
contributed to the results reported within this paper.


\addcontentsline{toc}{section}{References}
\bibliographystyle{plainurl}
\bibliography{nse-lpv-sdre}

\end{document}